\def\col{\hbox{col}}
\def\Rmin{\mathbb{R}_{\min}}
\def\R{\mathbb{R}}
\def\uc{\mathcal}

\documentclass[10pt, a4paper, twocolumn]{article} 

\usepackage{amsmath,amsfonts,amssymb,amsthm,titling,url,array}

\newtheorem{theorem}{Theorem}

\newtheorem{proposition}[theorem]{Proposition}
\newtheorem{lemma}[theorem]{Lemma}
\newtheorem{example}[theorem]{Example}
\usepackage{subfigure}

\usepackage{algorithm}
\usepackage{algpseudocode}

\usepackage[english]{babel} 

\usepackage{microtype} 

\usepackage{amsmath,amsfonts,amsthm} 

\usepackage[svgnames]{xcolor} 

\usepackage[hang, small, labelfont=bf, up, textfont=it]{caption} 

\usepackage{booktabs} 

\usepackage{lastpage} 

\usepackage{graphicx} 

\usepackage{enumitem} 
\setlist{noitemsep} 

\usepackage{tikz}
\usetikzlibrary{arrows}

\usepackage{sectsty} 
\allsectionsfont{\usefont{OT1}{phv}{b}{n}} 

\usepackage{geometry} 

\usepackage[T1]{fontenc} 
\usepackage[utf8]{inputenc} 

\usepackage{XCharter} 

\usepackage{fancyhdr} 
\fancypagestyle{firstpage}{ 
	\fancyhf{}
}

\newcommand{\authorstyle}[1]{{\large\usefont{OT1}{phv}{b}{n}\color{DarkRed}#1}} 

\newcommand{\institution}[1]{{\footnotesize\usefont{OT1}{phv}{m}{sl}\color{Black}#1}} 

\usepackage{titling} 

\newcommand{\HorRule}{\color{DarkGoldenrod}\rule{\linewidth}{1pt}} 

\pretitle{
	\vspace{-30pt} 
	\HorRule\vspace{10pt} 
	\fontsize{32}{36}\usefont{OT1}{phv}{b}{n}\selectfont 
	\color{DarkRed} 
}

\posttitle{\par\vskip 15pt} 

\preauthor{} 

\postauthor{ 
	\vspace{10pt} 
	\par\HorRule 
	\vspace{20pt} 
}

\usepackage{lettrine} 
\usepackage{fix-cm}	

\newcommand{\initial}[1]{ 
	\lettrine[lines=3,findent=4pt,nindent=0pt]{
		\color{DarkGoldenrod}
		{#1}
	}{}%
}

\usepackage{xstring} 

\newcommand{\lettrineabstract}[1]{
	\StrLeft{#1}{1}[\firstletter] 
	\initial{\firstletter}\textbf{\StrGobbleLeft{#1}{1}} 
}

\usepackage[backend=bibtex,style=numeric,natbib=true]{biblatex} 

\usepackage[autostyle=true]{csquotes} 

\title{Min-plus algebraic low rank matrix approximation: a new method for revealing structure in networks}

\author{\authorstyle{James Hook} 
	\newline\newline 
	\institution{University of Bath, United Kingdom}} 

\date{\today} 

\begin{document}

\maketitle 

\thispagestyle{firstpage} 


\lettrineabstract{In this paper we introduce min-plus low rank matrix approximation. By using min and plus rather than plus and times as the basic operations in the matrix multiplication; min-plus low rank matrix approximation is able to detect characteristically different structures than classical low rank approximation techniques such as Principal Component Analysis (PCA). We also show how min-plus matrix algebra can be interpreted in terms of shortest paths through graphs, and consequently how min-plus low rank matrix approximation is able to find and express the predominant structure of a network.}


\section*{Introduction}

Classical low rank matrix approximations, including techniques such as PCA, form the basis of many of the most commonly used algorithms in data science. These techniques presuppose that the data in question has some predominant linear structure. The low rank approximation extracts this structure, which can then be visualized or used to determine certain linear relationships that the data (approximately) satisfies. In this paper we extend the technique of low rank matrix approximation to the min-plus semiring. 

Min-plus algebra is the study of equations that are structured around the binary operations `taking the minimum' and `plus'. More formally min-plus algebra concerns the min-plus semiring $\mathbb{R}_{\min}=[\mathbb{R}\cup \{\infty\},\oplus,\otimes]$, where
$$
a\oplus b=\min(a,b), \quad a\otimes b=a+b, \quad \forall~a,b\in\mathbb{R}_{\min}.
$$

As we will show, min-plus low rank matrix approximation is closely related to classical low rank matrix approximation. But because we are working with the min-plus semiring, instead of a standard algebra, such as the field of real numbers, the min-plus low rank approximation is able to detect and express certain structures in the data, that could not be detected by classical approaches. Like in the classical case, min-plus low rank matrix approximation presupposes that the data has a linear structure, only now `linear' means `min-plus linear'.

Tropical algebra is the field of mathematics concerned with any semiring whose `addition' operation is max or min. For example the max-plus semiring or the max-times semiring. Max-plus algebra has many applications in dynamical systems and scheduling \cite{but10,heid2006}. Karaev and Miettinen have presented two approaches to max-times low rank matrix approximation \cite{capricorn,cancer}. Note that the max-times and min-plus semi-rings are isomorphic via $\phi:\mathbb{R}_{\max~\times}\mapsto \mathbb{R}_{\min~+}$ with $\phi(x)=-\log(x)$, but that this transformation does not preserve any standard norm, so that approximation in max-times is different to approximation in min-plus. Karaev and Miettinen show that max-times approximation can provide a useful alternative to classical Non-Negative Matrix Factorization (NNMF). In this paper we hope to show that min-plus approximation provides a useful tool for analyzing network structure from the point of view of pairwise shortest path distances.

The remainder of this paper is organized as follows. In Section 1 we introduce some standard definitions and basic results for min-plus matrix algebra, for a more thorough introduction to min-plus algebra see \cite{but10} and the references therein. In Sections 1.1 and 1.2 we introduce our formulation of min-plus low rank matrix approximation. In Section 2 we describe algorithms for solving min-plus regression and low rank matrix factorization problems. In Section 3 we apply our min-plus low rank matrix factorization algorithm to a small example network taken from Ecology.

\section{Min-plus matrices}\label{graphs}

A min-plus matrix is simply an array of entries from $\Rmin$, so that ${A}\in\Rmin^{n \times m}$ is an $n\times m$ array of elements which are each either a real number or $\infty$. Min-plus matrix multiplication is defined in analogy to the conventional plus-times case. For ${A}\in\Rmin^{n \times m}$ and ${B}\in\Rmin^{m \times d}$, we have $A\otimes B\in\Rmin^{n\times d}$ with
$$
(A\otimes B)_{ij}=\bigoplus_{k=1}^{m}\big(a_{ik}\otimes b_{kj}\big)=\max_{k=1}^{n}\big(a_{ik}+b_{kj}\big).
$$
For $A\in\Rmin^{n\times n}$ the \emph{precedence graph} $\Gamma(A)$ is defined to be the weighted directed graph with vertices $V=\{v(1),\dots,v(n)\}$ and an edge from $v(i)$ to $v(j)$ with weight $a_{ij}$, whenever $a_{ij}\neq \infty$. We write $A^{\otimes \ell}$ to mean the min-plus product of $A$ with itself $\ell$ times. 
\begin{proposition}
$\big(A^{\otimes \ell}\big)_{ij}=$ the weight of the minimally weighted path of length $\ell$, through $\Gamma(A)$, from $v(i)$ to $v(j)$. 
\end{proposition}
For $A\in\Rmin^{n\times n}$ the \emph{Kleene star} is defined by
\begin{equation}\label{kleene}
A^{\star}=I\oplus A \oplus A^{\otimes 2} \oplus\dots.
\end{equation}
where $I\in\Rmin^{n\times n}$ is the \emph{min-plus identity matrix} with $I_{ii}=0$ for $i=1,\dots,n$ and $I_{ij}=\infty$ for $i\neq j$.
\begin{proposition}
If $\Gamma(A)$ contains no negatively weighted cycles, then $A^{\star}$ exists and $\big(A^{\star}\big)_{ij}=$ the weight of the minimally weighted path through $\Gamma(A)$, from $v(i)$ to $v(j)$. Otherwise, \eqref{kleene} does not converge. 
\end{proposition} 
A matrix $A\in\Rmin^{n\times n}$ is \emph{idempotent} if $A\otimes A=A$. In general a function is idempotent if it acts as the identity on its image. Idempotent matrices form an important subset of min-plus matrices, with many special properties \cite{idempotent}. 
\begin{proposition}
Let $\mathcal{G}$ be a directed, weighted with non-negative edge weights, graph with vertices $V=\{v(1),\dots,v(n)\}$ and let $D\in\Rmin^{n\times n}$ with $d_{ij}=$ the weight of the minimally weighted path through $\mathcal{G}$ from $v(i)$ to $v(j)$. Then $D$ is idempotent.
\end{proposition} 
For ${A}\in\Rmin^{n \times m}$ the \emph{precedence bipartite graph} $\mathcal{B}(A)$ is the undirected, weighted, bipartite graph with vertices $X=\{x(1),\dots,x(n)\}$, $Y=\{y(1),\dots,y(m)\}$ and an undirected edge of weight $a_{ik}$ between $x(i)$ and $y(k)$, whenever $a_{ik}\neq\infty$.
\begin{proposition}
$\big(A\otimes A^{T}\big)_{ij}=$ the weight of the minimally weighted path through $\uc{B}(A)$ from $x(i)$ to $x(j)$.
\end{proposition}
For ${A}\in\Rmin^{n \times m}$ and ${B}\in\Rmin^{m \times d}$ the \emph{precedence tripartite graph} $\mathcal{T}(A,B)$ is the directed, weighted, tripartite graph with vertices $X=\{x(1),\dots,x(n)\}$, $Y=\{y(1),\dots,y(m)\}$ and $Z=\{z(1),\dots,z(d)\}$, an edge of weight $a_{ik}$ from $x(i)$ to $y(k)$, whenever $a_{ik}\neq\infty$ and an edge of weight $b_{kj}$ from $y(k)$ to $z(j)$, whenever $b_{kj}\neq\infty$.
\begin{proposition}
$\big(A\otimes B\big)_{ij}=$ the weight of the minimally weighted path through $\mathcal{T}(A,B)$ from $x(i)$ to $z(j)$.
\end{proposition}

\subsection{Min-plus low rank approximation of symmetric minimally weighted path matrices}

Let $\mathcal{G}$ be an undirected, weighted with non-negative edge weights, graph with vertices $V=\{v(1),\dots,v(n)\}$ and let $D\in\Rmin^{n\times n}$, with $d_{ij}=$ the weight of the minimally weighted path through $\mathcal{G}$ from $v(i)$ to $v(j)$. Now let $W=\{w(1),\dots,w(m)\}\subset V$ be a subset of the vertices of $\mathcal{G}$, which we call a set of \emph{waypoints}, then for $D_{W}=D(:,W)\in\Rmin^{n\times m}$, we have $\big(D_{W}\otimes D_{W}^{T}\big)_{ij}=$ the weight of the minimally weighted path through $\uc{B}(D_{W})$ from $x(i)$ to $x(j)$, or equivalently, $\big(D_{W}\otimes D_{W}^{T}\big)_{ij}=$ the weight of the minimally weighted path through $\uc{G}$, from $v(i)$ to $v(j)$, that goes via at least one waypoint in $W$. We define the rank-m \emph{actual waypoint approximation} of $D$ to be given by
\begin{equation}\label{actualw}
\min_{W\subset V,~|W|=m } \|D-D_{W}\otimes D_{W}^{T}\|_{F}.
\end{equation}
In practice, actual waypoint approximations tend to be quite inaccurate and the set of subsets of vertices is difficult to optimize over. Instead we use the rank-m \emph{virtual waypoint approximation} of $D$, which is defined by
\begin{equation}\label{virtualw}
\min_{A\in\Rmin^{n\times m}} \|D-A\otimes A^{T}\|_{F}.
\end{equation}
A solution to \eqref{virtualw} provides a bipartite graph $\uc{B}(A)$, with vertices $X=\{x(1),\dots,x(n)\}$ and $Y=\{y(1),\dots,y(m)\}$, 
 such that the weight of the minimally weighted path through $\uc{G}$, from $v(i)$ to $v(j)$, is approximately equal to the weight of the minimally weighted path through $\uc{B}(A)$ from $x(i)$ to $x(j)$, for all $i=1,\dots,n$, $j=1,\dots,m$.

\begin{example}\label{apfacegg}
Consider the following matrix $A\in\Rmin^{6\times 6}$. The precedence graph of $A$ is given in Figure~\ref{pregraph1}. We compute the shortest path distance matrix $D=A^{\star}$. 

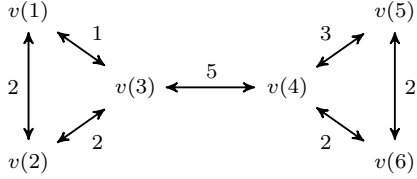
\begin{figure}
\centering
\begin{tikzpicture}[<->,>=stealth',shorten >=1pt,auto,node distance=1.5cm,main node/.style={font=\sffamily\footnotesize\bfseries}]

  \node[main node] (u1) {$v(1)$};
  \node[main node] (u2) at (0,-2) {$v(2)$};
    \node[main node] (u3) at (1.41,-1) {$v(3)$};

  \node[main node] (w1) at (3.41,-1) {$v(4)$};
 \node[main node] (w2) at (4.82,0) {$v(5)$};
 \node[main node] (w3) at (4.82,-2) {$v(6)$};

  \path[every node/.style={font=\sffamily\footnotesize}]

  (u1) edge [line width =0.25mm] node [above right]  {$1$} (u3)
 (u2) edge [line width =0.25mm] node [below right]  {$2$} (u3)    
 (u1) edge [line width =0.25mm] node [left] {$2$} (u2)
 
 (u3) edge [line width =0.25mm] node [above] {$5$} (w1)
 
 (w1) edge [line width =0.25mm] node [below left]  {$2$} (w3)
 (w2) edge [line width =0.25mm] node [right]  {$2$} (w3)    
 (w1) edge [line width =0.25mm] node [above left] {$3$} (w2);

 \end{tikzpicture}
 \caption{Precedence graph $\Gamma(A)$ of Example~\ref{apfacegg}.} \label{pregraph1}
 \end{figure}
\tiny
$$
A=\left[\begin{array}{cccccc} 0 & 2 & 1 & \infty & \infty & \infty \\ \cdot & 0 & 2 &  \infty & \infty & \infty \\ \cdot & \cdot & 0 & 5 & \infty & \infty \\ \cdot & \cdot & \cdot & 0 & 3 & 2 \\ \cdot & \cdot & \cdot & \cdot & 0 & 2 \\ \cdot & \cdot & \cdot & \cdot & \cdot & 0 \end{array}\right], \quad D=\left[\begin{array}{cccccc} 0 & 2 & 1 & 6 & 9 & 8 \\ \cdot & 0 & 2 &  7 & 10 & 9 \\ \cdot & \cdot & 0 & 5 & 8 & 7 \\ \cdot & \cdot & \cdot & 0 & 3 & 2 \\ \cdot & \cdot & \cdot & \cdot & 0 & 2 \\ \cdot & \cdot & \cdot & \cdot & \cdot & 0 \end{array}\right].
$$
    \normalsize
The waypoint set $W=\{v(3),v(4)\}$ yields the actual waypoint approximate factorization 

\tiny
$$
\left[\begin{array}{cc} 1 & 6 \\ 2 & 7 \\ 0 & 5 \\ 5 & 0 \\ 8 & 3 \\ 7 & 2 \end{array}\right]\otimes \left[\begin{array}{cccccc} 1 & 2 & 0 & 5 & 8 & 7 \\ 6 & 7 & 5 & 0 & 3 & 2 \end{array}\right]=\left[\begin{array}{cccccc} 2 & 3 & 1 & 6 & 9 & 8 \\ \cdot & 4 & 2 &  7 & 10 & 9 \\ \cdot & \cdot & 0 & 5 & 8 & 7 \\ \cdot & \cdot & \cdot & 0 & 3 & 2 \\ \cdot & \cdot & \cdot & \cdot & 6 & 5 \\ \cdot & \cdot & \cdot & \cdot & \cdot & 4 \end{array}\right],
$$
    \normalsize
 which results in a residual with Frobenius norm $9.5917$. Using Algorithm~\ref{newtomsymfac} we compute a virtual waypoint approximate factorization $F\otimes F^{T}\approx D$ with 
 \footnotesize
$$
F=\left[\begin{array}{cccccc} 0.4722  &  0.9722 &   0.2222  &  5.4444  &  9.0278 &   8.0278 \\
    7.7778   & 8.9778  &  6.7778 &   0.8889 &   1.0222  &  0.4222\end{array}\right],
    $$
    \normalsize
which results in a residual with Frobenius norm $4.5680$.

\end{example}

\subsection{Min-plus low rank approximation of general min-plus matrices}

The virtual waypoint approximation can be applied to general min-plus matrices. For $M\in\Rmin^{n\times d}$, the rank-m \emph{virtual waypoint approximation} of $M$ is defined by
\begin{equation}\label{virtualwAB}
\min_{A\in\Rmin^{n\times m},~B\in\Rmin^{m\times d}} \|M-A\otimes B\|_{F}.
\end{equation}
A solution to \eqref{virtualwAB} provides a tripartite graph $\uc{T}(A,B)$, with vertices $X=\{x(1),\dots,x(n)\}$, $Y=\{y(1),\dots,y(m)\}$ and $Z=\{z(1),\dots,z(d)\}$, such that the minimally weighted path through $\uc{T}(A,B)$ from $x(i)$ to $z(j)$ is approximately equal to $m_{ij}$, for all $i,j=1,\dots,n$.

\section{Algorithms for min-plus regression and low rank approximation}\label{algorithmsec}

An important prerequisite to matrix factorization is linear regression. For $A\in\Rmin^{n\times d}$ and $y\in\Rmin^n$ we seek
\begin{equation}\label{pregression}
\min_{x\in\Rmin^d}\|A\otimes x-y\|_{p},
\end{equation}
for some $p\in[1,\infty]$.

A function $f:\Rmin^{n}\mapsto\Rmin$ is \emph{min-plus convex} if for all $x,y\in\mathbb{R}_{\min}^{n}$ and $\lambda,\mu\in\mathbb{R}_{\min}$ such that $\lambda\oplus \mu=0$  we have
$$
f(\lambda\otimes x\oplus \lambda\otimes y)\leq \lambda\otimes f(x)\oplus \mu\otimes f(y).
$$
\begin{theorem}For $A\in\Rmin^{n\times d}$ and $y\in\Rmin^n$, the residual $r(x)=\|A\otimes x-y\|_{\infty}$ is min-plus convex.
\end{theorem}

\begin{theorem}
For $A\in\Rmin^{n\times n}$ and $y\in\Rmin^n$ let 
$$
\hat{x}=-\big(A^{T}\otimes(-y)\big), \quad x^{\ast}=\hat{x}\otimes \alpha
$$
where $\alpha=\max_{i}(y_{i}-(A\otimes \hat{x})_{i})/2$. Then
$$
x^{\ast}=\inf_{\leq}\big(\arg\min_{x\in\Rmin^d}|A\otimes x-y|_{\infty}\big)
$$
That is the infimum, with respect to the standard partial order $\leq$, of the set of optimal solutions.
\end{theorem}

\begin{example}\label{regegg}
Consider
$$
A=\left[\begin{array}{cc} 0 & 0 \\ 1 & 0 \\ 0 & 1 \end{array}\right], \quad y=\left[\begin{array}{cc}  0 \\ 1 \\ 1 \end{array}\right].
$$
\begin{figure*}
\begin{center}
\subfigure{\includegraphics[scale=0.38]{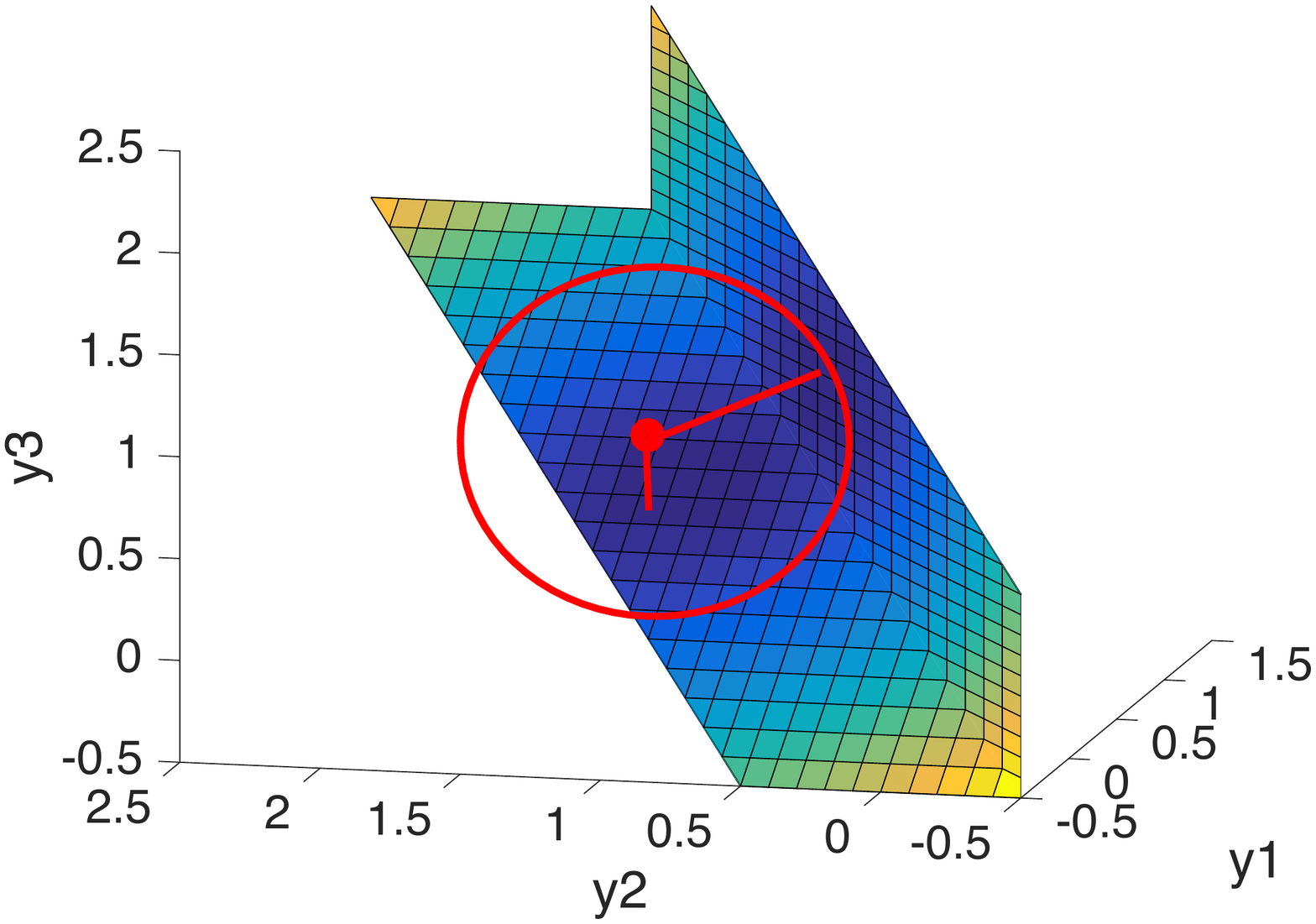}}  \hspace{0mm} \subfigure{\includegraphics[scale=0.38]{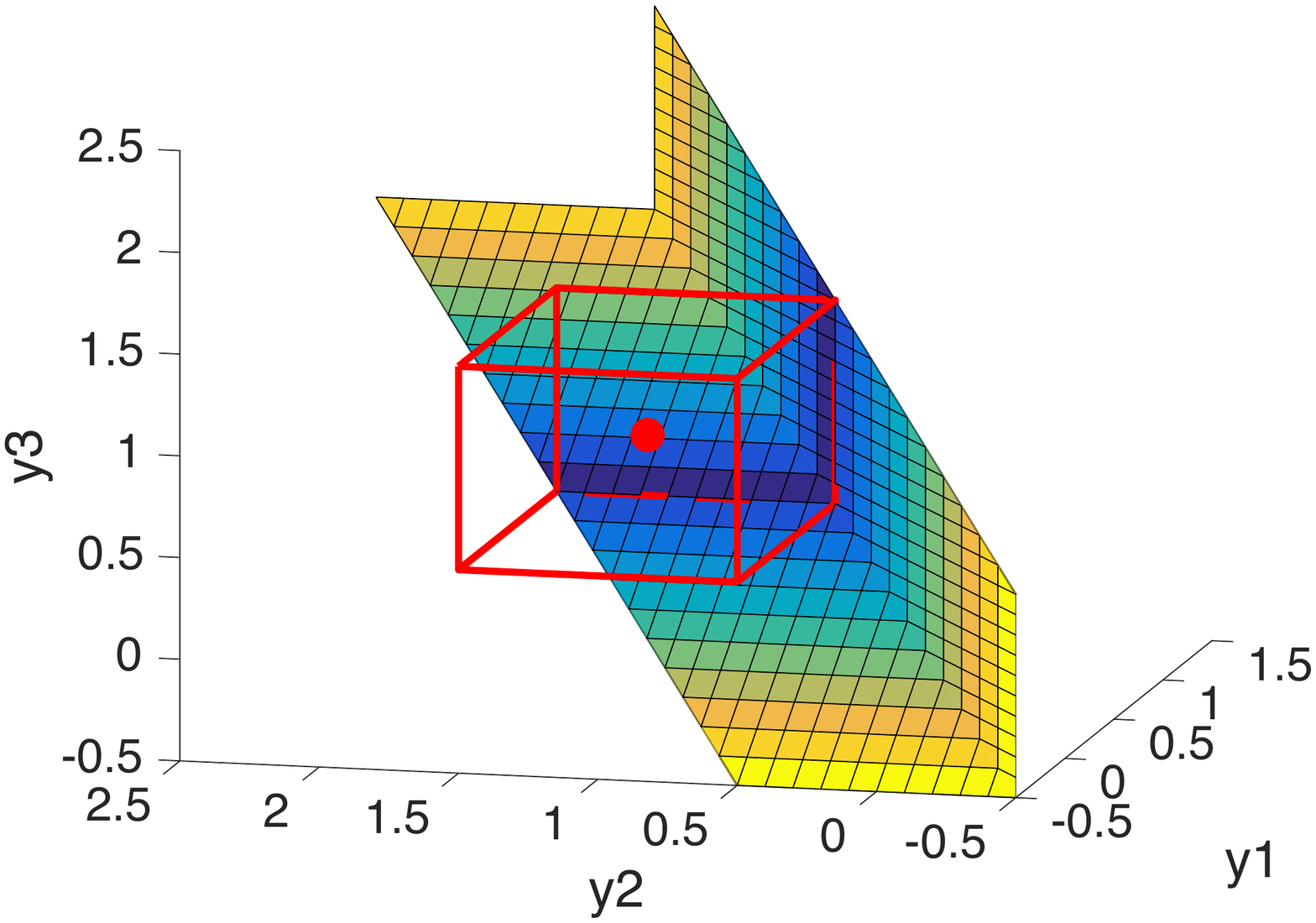}}
\caption{Column space geometry for the regression problem of Example~\ref{regegg}. Left $p=2$, right $p=\infty$.}	\label{colgeo1} 
\end{center}
\end{figure*}
Figure~\ref{colgeo1} displays the point $y$ and the column space $\col(A)=\{A\otimes x~:~x\in\Rmin^2\}$.  The solution to \eqref{pregression} is simply the closest point in $\col(A)$ to $y$, measured in the $p$-norm. For $p=2$ there are multiple local minima. For this example both local minima are global minima but typically this will not be the case. For $p=\infty$ there is a continuum of local minima. This set of minima is not convex with respect to plus and times but it is min-plus convex.

\end{example}

For $A\in\Rmin^{n\times d}$, $y\in\Rmin^n$ and $p=2$, the squared residual surface $r(x)^2=\|A\otimes x-y\|_{2}^2 $ is piecewise quadratic, continuous but non-differentiable. For $x\in\Rmin^d$ let
\begin{equation}
J(i,x)=\arg\min_{j=1}^{d} (a_{ij}+x_j ),
\end{equation}
and let $\tilde{J}(i,x)=\inf J(i,x)$, for $i=1,\dots,n$. Then we have
\begin{equation}
r(x)^2=p_{x}(x),
\end{equation}
where $p_{x}:\Rmin^{d}\mapsto\R$ is given by
\begin{equation}
p(x')=\sum_{i=1}^{n}\big(a_{i\tilde{J}(i,x)}+x'_{\tilde{J}(i,x)}-y_i \big)^2.
\end{equation}
Newton's method iteratively finds the minimum to the local quadratic piece
\begin{equation}
\mathcal{N}(x)=\min_{x'\in\Rmin^{d}}p_{x}(x'),
\end{equation}
which is given by
\begin{equation}
\mathcal{N}(x)_{k}={\big(\sum_{i:\tilde{J}(i,x)=k}y_i-a_{ik}\big)}~\Big/{\big(\sum_{i:\tilde{J}(i,x)=k}1\big)}.
\end{equation}
However, since the residual surface is non-differentiable Newton's method isn't guaranteed to converge to a local minimum. Therefore we propose using a Newton directed line search. In order to search efficiently it is important to take into account the discontinuities in the residuals derivative. If on one iteration the line search minimum $x(k)$ is found to lie on the discontinuity surface $\Sigma$, then the following Newton's step must be restricted to directions tangental to $\Sigma$ at $x(k)$. See Algorithm~\ref{NDLS}.

\begin{algorithm}
\caption{ \label{NDLS}
Given $A\in\Rmin^{n\times n}$, $y\in\Rmin^n$ and an initial guess $x(0)\in\Rmin^d$, this algorithm returns a local minimum of $r(x)=\|A\otimes x-y\|_{2}$.}
\begin{algorithmic}[1]
\While{ not converged}
\State compute $\mathcal{N}\big(x(k)\big)$, restricting to $T_{x(k)}\Sigma$ if $x(k)\in\Sigma$
\State find $\lambda_{\min}=\arg\min_{\lambda}r\Big(L\big(x(k),\lambda\big)\Big)$, where $L\big(x(k),\lambda\big)=\lambda\mathcal{N}\big(x(k)\big)+(1-\lambda)x(k)$
\If{$\lambda_{\min}=1$}
\State return $x=\mathcal{N}\big(x(k),\sigma(k)\big)$
\EndIf
\State $x(k+1)=L\big(x(k),\lambda_{\min}\big)$
\EndWhile
\end{algorithmic}
\end{algorithm}

The Newton update (line 2) can be computed with cost $\uc{O}(nd)$ per iteration. The line search (line 3) can be computed with cost $\uc{O}\big(nd\log(nd)\big)$ per iteration, as follows. The $i$th component of the product $A\otimes L\big(x(k),\lambda\big)$ is a piecewise affine function of $\lambda$, with up to $d$ points of non-differentiability, for $i=1,\dots,n$. Thus $r\Big(L\big(x(k),\lambda\big)\Big)^2$ is piecewise quadratic, with at most $nd$ points of non-differentiability. We begin by finding the points of non-differentiability, then sort them, then carry out the line search through the quadratic pieces. If the minimum is attained at a non-differentiability point then we know that the line search minimum is attained on the discontinuity surface $\Sigma$. The cost $\uc{O}\big(nd\log(nd)\big)$ is the worst case cost, associated with sorting the maximum possible number of non-differentiability points.  The infimum of the set of optimal solutions to the $p=\infty$ problem provides a good choice for the initial condition. 

Algorithm~\ref{NDLS} enables a simple alternating method for computing a non-symmetric approximate factorizations of a general min-plus matrix $M\in\Rmin^{n\times d}$. The only difficulty is choosing an initial factorization to work from. One possibility is to run the kmeans clustering algorithm to find $m$ centers that approximate the $d$ columns of $M$. We use these centers as the initial LHS factor. We then use the formula for the infimum solution of the $p=\infty$ regression problem to fit an initial RHS factor. See Algorithm~\ref{nonsymfac}.

\begin{algorithm}
\caption{ \label{nonsymfac}
Given $M\in\Rmin^{n\times d}$ and $0<m\leq \min(n,d)$, this algorithm returns an approximate factorization $A\otimes B\approx M$ with $A\in\Rmin^{n\times m}$, $B\in\Rmin^{m\times d}$.}
\begin{algorithmic}[1]
\State compute $A(0)=\hbox{kmeans}(M,k)$
\State set $B(0)_{\cdot j}=\inf \arg \min_{x\in\Rmin^{d}}\|A(k-1)\otimes x-M_{\cdot j}\|_{\infty}$, for $j=1,\dots,m$.
\While{ not converged}
\State set $B(k)_{\cdot j}=\hbox{ALG1}\big(A(k-1),M_{\cdot j},B(k-1)_{\cdot j}\big)$, for $j=1,\dots,m$.
\State set $A(k)_{i\cdot}=\hbox{ALG1}\big(B(k-1)^{T},M_{i\cdot}^{T},A(k-1)_{i\cdot}^{T}\big)^{T}$, for $i=1,\dots,n$.
\If {$[A(k),B(k)]=[A(k-1),B(k-1)]$}
\State return $[A,B]=[A(k),B(k)]$
\EndIf
\EndWhile
\end{algorithmic}
\end{algorithm}

Fitting the initial LHS factor (line 1) with kmeans has cost $\uc{O}(nmk)$ per kmeans iteration. Fitting the initial RHS factor (line 2) has cost $\uc{O}(nmk)$. We call Algorithm~\ref{NDLS} to update each column in the RHS factor (line 4), using the previous value as the initial guess. This has cost $\uc{O}\big(mnk\log(nk)\big)$ per iteration. Similarly updating the LHS factor (line 5) has cost $\uc{O}\big(nmk\log(mk)\big)$ per iteration.

Computing a symmetric approximate factorization for a symmetric shortest path distance matrix is a little more difficult. We cannot use Algorithm~\ref{nonsymfac} as we do not have any compatible way of enforcing symmetry in the factors at each step. Instead we apply Newton's method, using the whole of the approximate factor as the iterate. For $D\in\Rmin^{n\times n}$ the squared residual surface $r(F)^2=\|F\otimes F^{T}-D\|_{F}^{2}$
is piecewise quadratic, continuous but non-differentiable. For $F\in\Rmin^{n\times m}$ let
$$
K(i,j,F)=\arg\min_{k=1}^{d}f_{ik}+f_{jk},
$$
and let $\tilde{K}(i,j,F)=\inf K(i,j,F)$, for $i,j=1,\dots,n$. Then we have
\begin{equation}\label{sfresidual3}
r(F)^2=q_{F}(F),
\end{equation}
where $q_{F}:\Rmin^{n\times m}\mapsto \R$ is given by
\begin{equation}\label{sfresidual2}
q_{F}(F')=\sum_{ij=1}^n\big(f'_{i\tilde{K}(i,j,F)}+f'_{j\tilde{K}(i,j,F)}-d_{ij}\big)^2.
\end{equation}
As in the case of the regression problem, Newton's method finds the minimum to the local quadratic piece
\begin{equation}\label{nretonsym}
\uc{N}(F)=\arg\min_{F'\in\Rmin^{n\times m}}q_{F}(F').
\end{equation}
Define $\uc{J}_{F}:\Rmin^{n\times m}\mapsto\Rmin^{n\times m}$ by
\begin{equation}\label{jacobi}
\uc{J}_{F}(F')_{ik}=\frac{d_{ii}{\bf 1}_{ik}+\sum_{j\neq i : \tilde{K}(i,j,F)=k}d_{ij}-f'_{jk}}{2{\bf 1}_{ik}+\sum_{j\neq i : \tilde{K}(i,j,F)=k}1},
\end{equation}
where ${\bf 1}_{ik}=1$ if $\tilde{K}(i,i,F)=k$ and ${\bf 1}_{ik}=0$ otherwise. The map $\uc{J}_{F}$ is simply the result of applying one iteration of Jacobi's method to the normal equations associated to the linear least squares formulation of \eqref{nretonsym}.
\begin{lemma}
Let $F(0)=F$ and let $F(t+1)=\uc{J}_{F}\big(F(t)\big)$ for $t=0,1,\dots$ then
$$
\lim_{t\rightarrow\infty}F(t)=\uc{N}(F).
$$
\end{lemma}
Therefore we can compute Newton's method updates iteratively using $\uc{J}$. However, as in the case of the regression problem, Newton's method is not guaranteed to converge. One possibility is to use a Newton directed line search, as we did in Algorithm~\ref{NDLS}. However, this would require us to store and sort $\uc{O}(n^2m)$ breakpoints, which presents a huge memory requirement even for modestly large matrices.  Instead we propose using approximate Newton updates with undershooting. By using a small fixed number of $\uc{J}$ iterations we can cheaply approximate the Newton step. We then update by moving to a point somewhere between the previous state and the result of our approximate Newton computation. By gradually reducing the length of the step we can avoid getting stuck in the periodic orbits that prevent standard Newton's method from converging. As in the non-symmetric case the choice of initial factorization is very important. One possibility is to take an actual waypoint factorization, as defined in \eqref{actualw}, using a randomly chosen subset of $m$ vertices as the waypoints. See Algorithm~\ref{newtomsymfac}. Tuning the number of Jacobi iterations at each step $t$, as well as the stepping parameter $\mu$ and the convergence criteria are important factors for the algorithms performance, which we hope to explore in detail in future work.

\begin{algorithm}
\caption{ \label{newtomsymfac}
Given a symmetric shortest path distance matrix $D\in\Rmin^{n\times n}$ and $0<k\leq d$, this algorithm returns an approximate factorization $F\otimes F^{T}\approx D$ with $F\in\Rmin^{n\times k}$. Parameters are the number of Jacobi iterations per step $t\in\mathbb{N}$ and the shooting factor $\mu\in\mathbb{R}_{+}$. These parameters may be allowed to vary during the copmutation. }
\begin{algorithmic}[1]
\State draw uniformly at random $\{w_{1},\dots,w_{m}\}\subset\{1,\dots,n\}$
\State set $F(0)=D_{W}$
\While{ not converged}
\State compute $\hat{\uc{N}}\big(F(k-1)\big)=\uc{J}^{t}_{F(k-1)}\big(F(k-1)\big)$
\State set $F(k)=\mu\hat{\uc{N}}\big(F(k-1)\big)+(1-\mu)F(k-1)$
\EndWhile
\State return $F$
\end{algorithmic}
\end{algorithm}

Formulating the map $\uc{J}_{F}$ has cost $\uc{O}(n^2m)$ and applying it has cost $\uc{O}(n^2)$. Thus the approximate Newton computation (line 4) has cost $\uc{O}\big(n^2(m+t)\big)$, where $t$ is the number of Jacobi iterations used at each step.

\section{Example: Latent factor analysis of dolphin social network}\label{numerics}

In this example we examine a small social network to illustrate min-plus low rank matrix approximation's ability to extract and visualize predominant structure in network data. We use the dolphin social network presented in \cite{dolphins}. This network consists of 62 vertices, each of which represents a different dolphin, with an edge connecting two dolphins if they are observed to regularly interact. This small social network is frequently used to test or illustrate data analysis techniques.  We construct the adjacency matrix $A\in\{0,1\}^{62\times 62}$, with $a_{ij}=1$, if and only if dolphin $i$ and dolphin $j$ are connected. We then compute the distance matrix $D\in\Rmin^{62\times 62}$, with $d_{ij}=$ the length of the shortest path through the network from dolphin $i$ to $j$.

\begin{figure*}[h!]
\begin{center}
\subfigure{\includegraphics[scale=0.35]{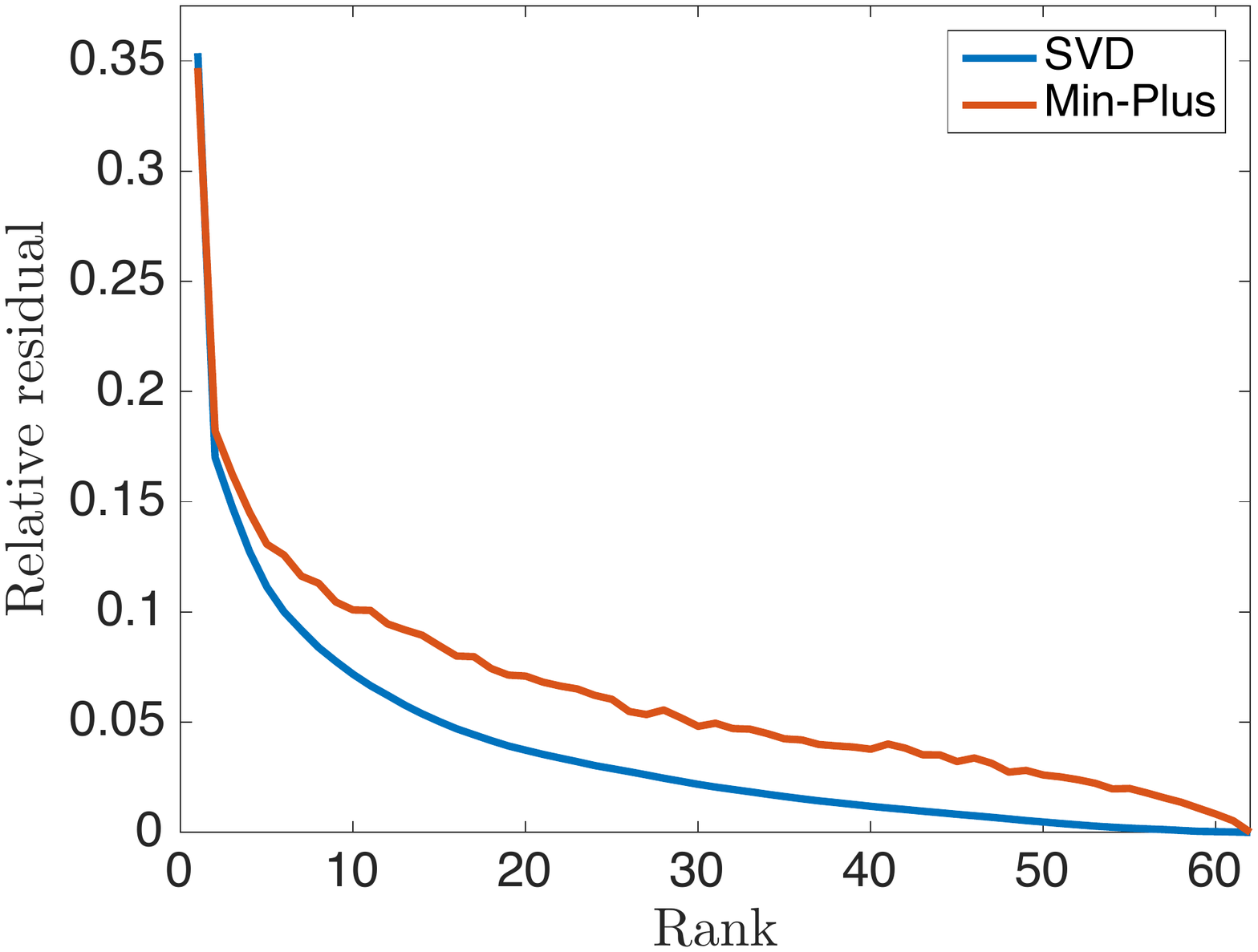}}  \hspace{15mm} \subfigure{\includegraphics[scale=0.35]{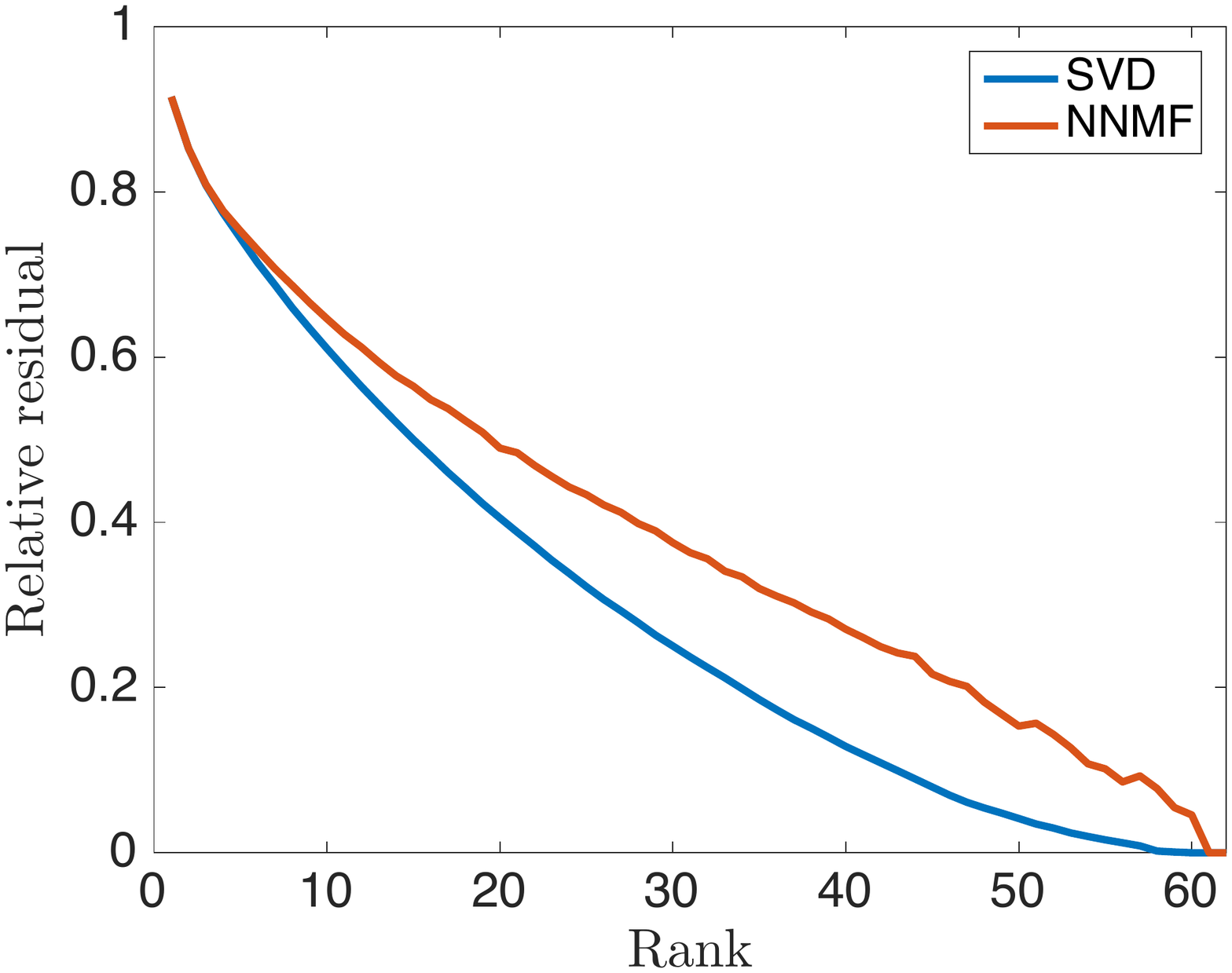}}
\caption{Comparison of residual errors for approximate factorization. Left: $D$, Right: $A$.}\label{dolphins0}
\end{center}
\end{figure*}

For $d=1,\dots,n$ we compute a rank-$d$ min-plus low rank matrix approximation of $D$ using Algorithm~\ref{nonsymfac}. We use the parameters $t=5$, $\mu=0.5$ and stop the algorithm after 100 iterations. We repeatedly run the algorithm 100 times with different initial conditions, then save the factorization that has the smallest residual error. Figure~\ref{dolphins0} displays a plot of the relative residuals for the min-plus factorization, given by $\|D-F\otimes F^{T}\|_{F}/\|D\|_{F}$, as well as for the truncated SVD as a comparison. In this example the truncated SVD has a smaller error for intermediate values of the rank, but for very small or very large ranks the residuals are nearly identical. It is not clear yet, whether the better performance of SVD for intermediate values is due to the data having a closer classical linear structure or Algorithm~\ref{nonsymfac} falling to find a good minimizer. 

A major advantage of the min-plus low-rank factorization over truncated SVD is the interpretability of the factors. We can think of the columns of $F$ as representing neighborhoods or waypoints. If $f_{ik}$ is small then dolphin $i$ is close to neighborhood $k$, and therefore dolphin $i$ will be close to any other dolphin that is also close to neighborhood $k$. Otherwise if $f_{ik}$ is large then dolphin $i$ is far from neighborhood $k$ and will not be close to any dolphin that is close to neighborhood $k$, unless they share some other mutually close neighborhood. 

\begin{figure*}[t]
\begin{center}
\subfigure{\includegraphics[scale=0.45]{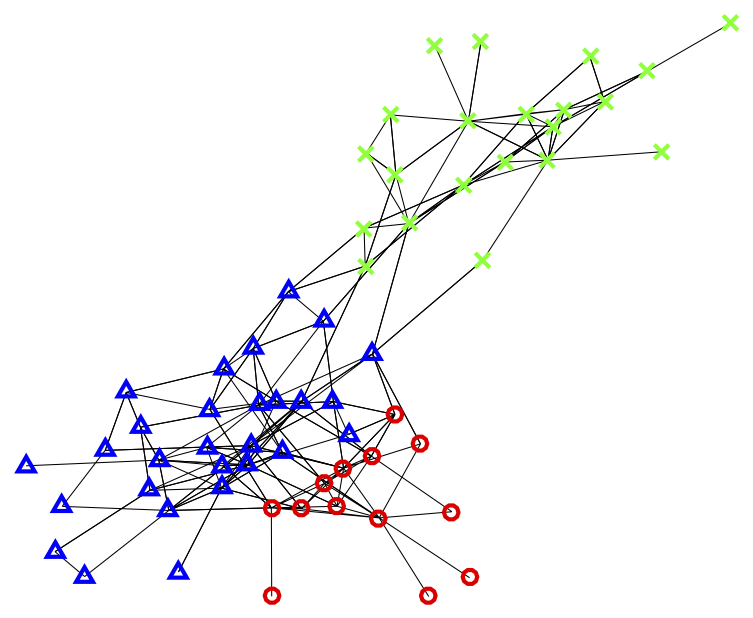}}  \hspace{15mm}  \subfigure{\includegraphics[scale=0.45]{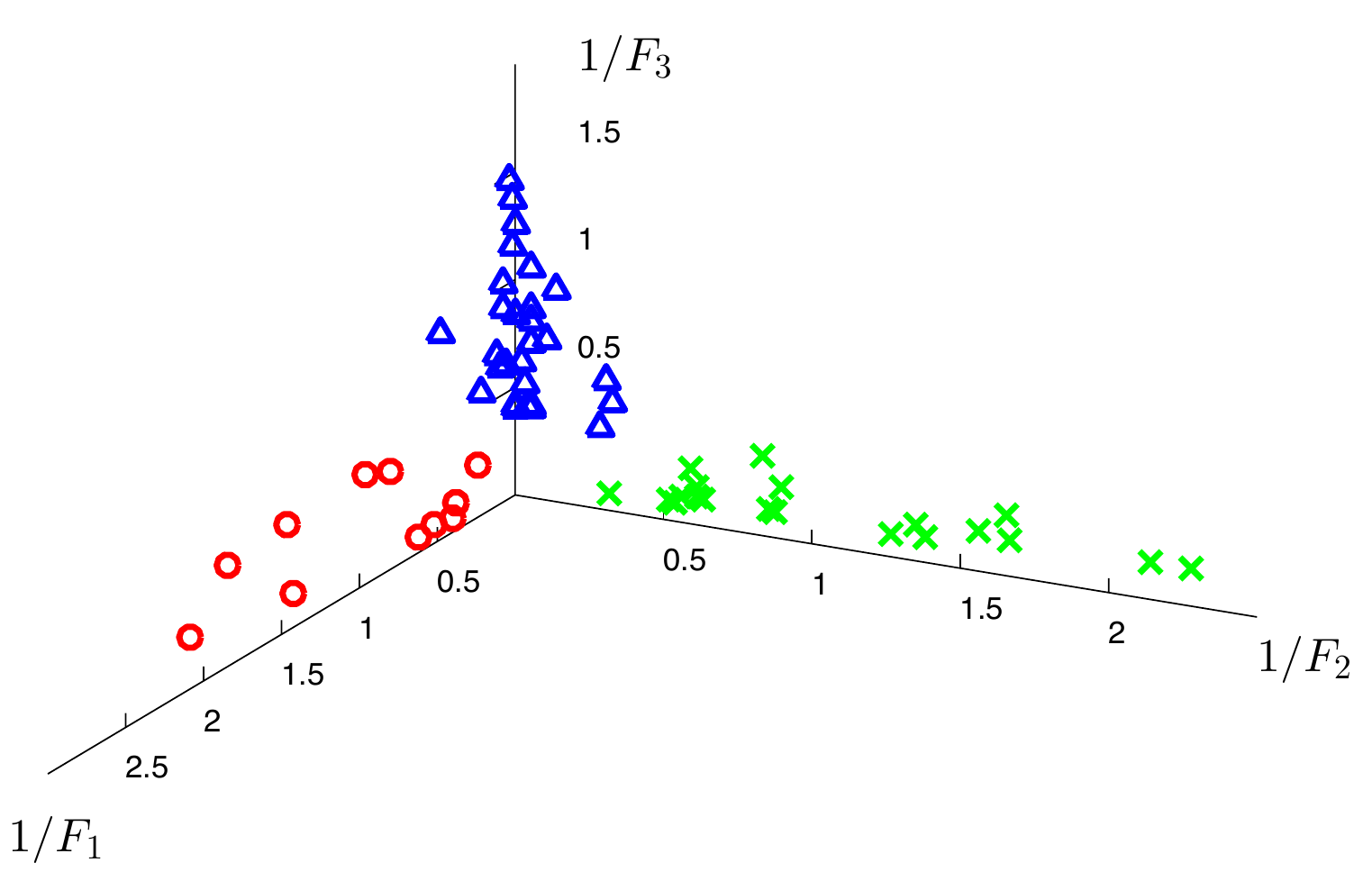}} 
\caption{Left: Dolphin social graph, Right: Min-plus latent factors.}\label{dolphins1}
\end{center}
\end{figure*}

Just as in PCA, the rows of $F$ can be thought of as latent factors that parametrize the rows of $D$. Equivalently, the $i$th row $F_{i\cdot}$ encapsulates information about dolphin $i$'s position in the network, so we can study the structure of the network by examining $\{F_{i\cdot}~:~i=1,\dots,n\}$, which is simply a scattering of points in $\R^3$. Figure~\ref{dolphins1} displays the dolphin social graph as well as the rows of $F$. Note that we have plotted the reciprocals of the entries in $F$, so that a large value of $1/f_{ik}$ indicates that dolphin $i$ is close to neighborhood $k$. The dolphins have then been color coded according to their closest neighborhood. Comparing the network to the scattering of points, it is clear that the min-plus factorization has captured the predominant structure of the graph.

\subsubsection{Comparison with non-negative matrix factorization}

Non-negative matrix factorization (NNMF) is a popular technique for community detection, which can be interpreted as the maximum a posteriori estimate for the inverse problem of inferring a stochastic block model to explain the network structure. Like min-plus low rank approximation, a low rank approximate non-negative matrix factorization tends to result in a larger residual that a truncated SVD of the same rank but is preferred in some situations becuase of the better interpretability of the factors. Figure~\ref{dolphins0} displays a plot of the relative residuals for the truncated SVD and non-negative matrix factorizations of $A$ as a function of rank. A fundamental difference between the NNMF approach and our min-plus low rank matrix approximation, is that the NNMF is applied directly to the adjacency matrix $A$, whilst our approach is to factorize the distance matrix $D$. This means that our method is sensitive to the indirect connections between vertices, which NNMF is oblivious to. In some cases these connections may not be of interest, in which case NNMF is a good choice. However, in many applications such connections are extremely important. Suppose for example that a message or contagion was spread through the network, then the distances between non-directly connected dolphins would need to be taken into account.

\section{Conclusion}

We have introduced min-plus low rank matrix approximation. The small example in Section~\ref{numerics} demonstrates that min-plus low rank matrix approximation is able to detect and express predominant networks structure in a novel way that could be useful in a range of fields. In further work we hope to develop some useful techniques based on min-plus low rank approximation aimed at specific network analysis applications.


\printbibliography[title={Bibliography}] 

\end{document}